\documentclass[a4paper,12pt]{article}
\usepackage[T1]{fontenc}
\usepackage{mathptmx,courier,pifont}
\usepackage{amsmath,amssymb,amsfonts}        % moduli per le estensioni ed i fonts AMS
\usepackage[dvips]{graphicx}
\usepackage{psfrag}
\usepackage{latexsym}
\usepackage{color}
\usepackage{url}
\def\RR{\hbox{I\kern-.2em\hbox{R}}}
\def\NN{\hbox{I\kern-.2em\hbox{N}}}
\def\ds{\displaystyle}

\def\E{\mbox{E}}

\newcommand{\eqnsection}{
   \renewcommand{\theequation}{{\thesection.\arabic{equation}}}
   \makeatletter
   \csname @addtoreset\endcsname{equation}{section}
   \makeatother}
\eqnsection

\title{Finite difference schemes on quasi-uniform grids for BVPs on infinite intervals}

\author{Riccardo Fazio\thanks{Corresponding author home-page: http://mat521.unime.it/$\sim$fazio} \ \ and Alessandra Jannelli \\
Department of Mathematics and Computer Science \\
University of Messina \\
Viale F. Stagno D'Alcontres 31, 98166 Messina, Italy \\
rfazio@unime.it\ \ \ \ \ ajannelli@unime.it}
\date{Submitted to Journal Computational and Applied Mathematics on May 7, 2013 and in revised form November 27, 2013 and \today}
\begin{document}
\maketitle

\begin{abstract}
The classical numerical treatment of boundary value problems defined on infinite intervals is to replace the boundary conditions at infinity by suitable boundary conditions at a finite point, the so-called truncated boundary.
A truncated boundary allowing for a satisfactory accuracy of the numerical solution has to be determined by trial and errors and this seems to be the weakest point of the classical approach.
On the other hand, the free boundary approach overcomes the need for a priori definition of the truncated boundary. 
In fact, in a free boundary formulation the unknown free boundary can be identified with a truncated boundary and being unknown it has to be found as part of the solution. 

In this paper we consider a different way to overcome the introduction of a truncated boundary, namely non-standard finite difference schemes defined on quasi-uniform grids.
A quasi-uniform grid allows us to describe the infinite domain by a finite number of intervals.
The last node of such grid is placed on infinity so that right boundary conditions
are taken into account exactly.
We apply the proposed approach to the Falkner-Skan model and to a problem of interest in foundation engineering.
The obtained numerical results are found in good agreement with those available in literature.
Moreover, we provide a simple way to improve the accuracy of the numerical results using Richardson's extrapolation.
Finally, we indicate a possible way to extend the proposed approach to boundary value problems defined on the whole real line.
\end{abstract}



%\pagebreak
\section{Introduction}\label{S:intro}
The classical numerical treatment of
boundary value problems (BVPs) on infinite intervals is to replace the original problem by one defined
on a finite interval, where a finite value, the so-called
truncated boundary, is used instead of infinity (see, for instance, Collatz \cite[pp. 150-151]{Collatz} or Fox \cite[p. 92]{Fox}). 
For the accepted numerical solution, the value of the truncated boundary is varied until the computed results stabilize, at least, to a prefixed number of significant digits. 
%For instance, in the case of the von Karman swirling flows, the opportunity to use values of $ x_{\infty} $ up to $ 200 $ was reported by Lentini and Keller \cite{LentiniK80} in order to investigate a fourth branch of the flow.
However, a truncated boundary allowing for a satisfactory accuracy
of the numerical solution has to be determined by trial and errors
and this seems to be the weakest point of the classical approach.
Hence, a priori definition of the truncated boundary
was indicated by Lentini and Keller \cite{Lentini:BVP:1980} as an important area of research.

A theory for defining asymptotic boundary conditions to be imposed
at the truncated boundary has been developed by de Hoog and Weiss
\cite{deHoog:1980:ATB}, Lentini and Keller \cite{Lentini:BVP:1980} and Markowich
\cite{Markowich:TAS:1982,Markowich:ABV:1983}. 
%See also the related work by: Markowich \cite{Mark83}, Markowich and Ringhofer \cite{MarRin}, Schmeiser \cite{Scheme} and Mattheij \cite{Mattheij}. 
%In the last decade t
The asymptotic boundary conditions have been applied successfully to the numerical approximation of the so-called \lq \lq connecting orbits\rq \rq \ problems of dynamical systems, see Beyn \cite{Beyn:1990:GBN,Beyn:1990:NCC,Beyn:1992:NMD}. %Champneys and Kuznetsov \cite{Champneys}, Doedel \cite{Doedel86}, Doedel and Kernevez \cite{Doedel89}, Friedman and Doedel \cite{Friedman91,Friedman93}, Schecter \cite{Schecter93,Schecter95} and Sandstede \cite{Sandstede97}; more on this topic in the last section).
Those problems are of interest, not only in connection with dynamical systems, but also in the study of traveling wave solutions of partial differential equations of parabolic and hyperbolic type as shown by Beyn \cite{Beyn:1990:NCC}, Friedman and Doedel \cite{Friedman:1991:NCC}, Bai et al. \cite{Bai:1993:NCH}, and Liu et al. \cite{Liu:1997:CCH}.

A free boundary formulation was proposed by Fazio
\cite{Fazio:1992:BPF} where the unknown free boundary was identified with a truncated boundary. 
In this approach the free boundary is unknown and has to be found as part of the solution. 
This free boundary approach overcomes the need for a priori definition of the truncated boundary. 
The free boundary formulation has been applied to: the Blasius
problem \cite{Fazio:1992:BPF}, the Falkner-Skan model \cite{Fazio:1994:FSEb}, a model describing the flow  of an incompressible fluid over a slender parabola of revolution \cite{Fazio:1996:NAN}, a connecting orbit problem \cite{Fazio:2002:SFB}, and a problem in foundation engineering \cite{Fazio:2003:FBA}.

%\textcolor{red}{Citare \cite{Grosch:NSP:1977}.}
A different way to avoid the definition of a truncated boundary is to apply coordinate transforms.
The idea of mapping an infinite geometry into a finite one is not original.
For example, van de Vooren and Dijkstra \cite{vandeVooren:1970:NSS} applied coordinate transformations to the numerical solution of laminar flow past a flat plate, Botta et al. \cite{Botta:1972:NSN} and Davis \cite{Davis:1972:NSN} applied similar techniques to laminar flow past a parabola.
Coordinate transforms have been applied to the numerical solution of ordinary and partial differential equations on unbounded domains, see Grosch and Orszag \cite{Grosch:NSP:1977}, Boyd \cite{Boyd:2001:CFS} or Koleva \cite{Koleva:NSH:2006}.
Here we consider finite difference schemes on quasi-uniform grids, defined by coordinate transforms, applied to the numerical solution of BVPs defined on infinite intervals.
The novelty of our approach is that we define non-standard finite differences for the original problem on the infinite domain, whereas Grosch and Orszag transform the governing model and apply the classical finite difference or shooting methods on the transformed finite domain.
In the following sections we consider two test problems.
The first is the Falkner-Skan model of boundary layer theory.
The last one is a problem of interest in foundation engineering.
The obtained numerical results are found in good agreement with those available in literature.
Moreover, we have applied Richardson's extrapolation in order to improve the accuracy of the numerical results.
Preliminary numerical results on the main topic of this paper were presented at the ENUMATH 2013 conference \cite{Fazio:2014:QUG}.

We notice that the free boundary approach or a quasi-uniform grid strategy are as simple as the classical truncated boundary one in contrast with the asymptotic boundary approach, in this context see also the opinion expressed by J. R. Ockendon \cite{Ockendon}.

In the last section, we point out some conclusions supported by the
evidences of the present work and indicate a possible way to
extend the proposed approach to BVPs defined on the whole real
line. 

\section{Finite differences on quasi-uniform grids}\label{S:quniform}
Let us consider the smooth strict monotone quasi-uniform maps $x = x(\xi)$, the so-called grid generating functions,
\begin{equation}\label{eq:qu1}
x = -c \cdot \ln (1-\xi) \ ,
\end{equation}
and
\begin{equation}\label{eq:qu2}
x = c \frac{\xi}{1-\xi} \ ,
\end{equation}
where $ \xi \in \left[0, 1\right] $, $ x \in \left[0, \infty\right] $, and $ c > 0 $ is a control parameter.
So that, a family of uniform grids $\xi_n = n/N$ defined on interval $[0, 1]$ generates one parameter family of quasi-uniform grids $x_n = x (\xi_n)$ on the interval $[0, \infty]$.
The two maps (\ref{eq:qu1}) and (\ref{eq:qu2}) are referred as logarithmic and algebraic map, respectively. 
As far as the authors knowledge is concerned, van de Vooren and Dijkstra \cite{vandeVooren:1970:NSS} were the first to use these kind of maps. 
We notice that more than half of the intervals are in the domain with length approximately equal to $c$ and 
$x_{N-1} = c \ln N$ for (\ref{eq:qu1}),
while $ x_{N-1} \approx c N $ for (\ref{eq:qu2}).
For both maps, the equivalent mesh in $x$ is nonuniform with the
most rapid variation occurring with $c \ll x$.
The logarithmic map (\ref{eq:qu1}) gives slightly better resolution near $x = 0$ than the
algebraic map (\ref{eq:qu2}), while the algebraic map gives much better resolution than the
logarithmic map as $x \rightarrow \infty$. 
In fact, it is easily verified that
\[
-c \cdot \ln (1-\xi) < c \frac{\xi}{1-\xi} \ ,
\]
for all $\xi$, see figure \ref{fig:m1N20} below.

The problem under consideration can be discretized by introducing a uniform grid $ \xi_n $ of $N+1$ nodes in $ \left[0, 1\right] $ with $\xi_0 = 0$ and $ \xi_{n+1} = \xi_n + h $ with $ h = 1/N $, so that $ x_n $ is a quasi-uniform grid in $ \left[0, \infty\right] $. 
The last interval in (\ref{eq:qu1}) and (\ref{eq:qu2}), 
namely $ \left[x_{N-1}, x_N\right] $, is infinite but the point $ x_{N-1/2} $ is finite, because the non integer nodes are defined by 
\[
x_{n+\alpha} = x\left(\xi=\frac{n+\alpha}{N}\right) \ ,
\]
with $ n \in \{0, 1, \dots, N-1\} $ and $ 0 < \alpha < 1 $.
%In this way we have defined also the ghost cell $ x_{-1} $.  
This maps allow us to describe the infinite domain by a finite number of intervals.
The last node of such grid is placed on infinity so right boundary condition
is taken into account correctly.

Figure \ref{fig:m1N20} shows the two quasi-uniform meshes $x=x_n$, $n = 0, 1, \dots , N$ defined by (\ref{eq:qu1}) and by (\ref{eq:qu2}) with $c=5$ and $N=20$.
\begin{figure}[!hbt]
\centering
\psfrag{I}[][]{$x_0$} 
\psfrag{xxx}[l][]{$x \rightarrow \infty$} 
\psfrag{P}[][]{$x_{19}$} 
\psfrag{A}[][]{} 
\framebox{\includegraphics[width=.9\textwidth]{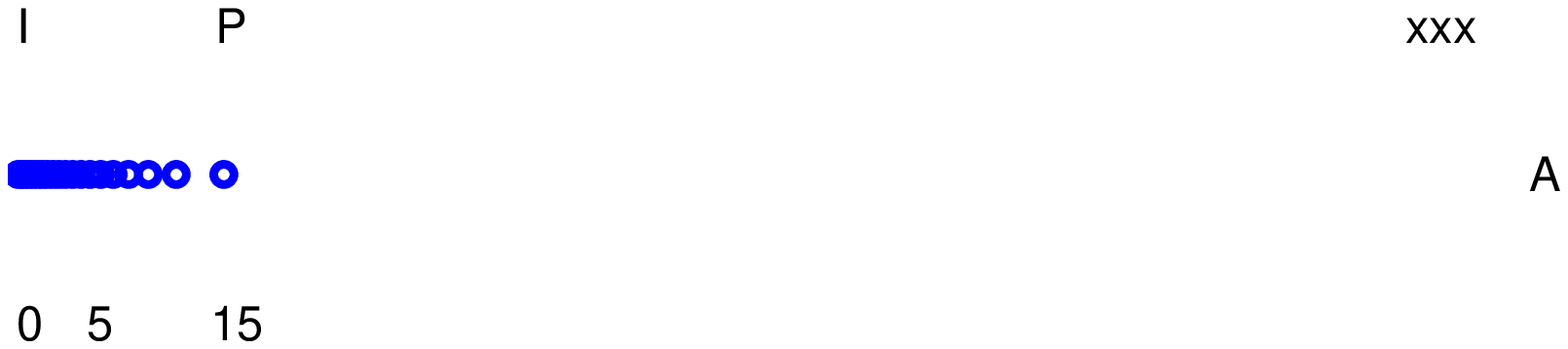}
} \\
\framebox{\includegraphics[width=.9\textwidth]{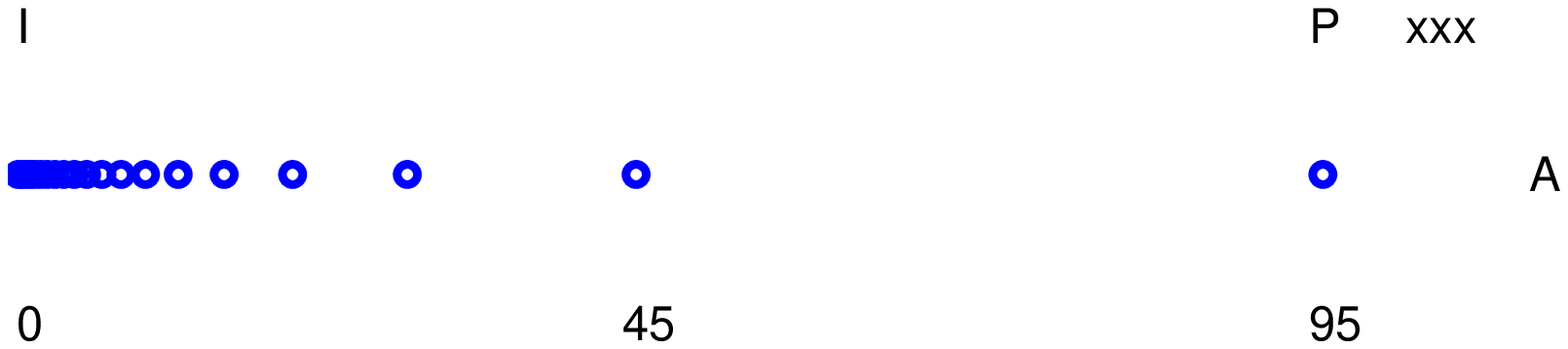}
}
\caption{\it Quasi-uniform meshes: top frame for (\ref{eq:qu1}) and bottom frame for (\ref{eq:qu2}). We notice that, in both cases, the last mesh-point is $x_N = \infty$.}
\label{fig:m1N20}
\end{figure}

We can define the values of $u(x)$ on the mid-points of the grid
\begin{equation}
u_{n+1/2} \approx \frac{x_{n+1}-x_{n+1/2}}{x_{n+1}-x_n} u_n + \frac{x_{n+1/2}-x_{n}}{x_{n+1}-x_n} u_{n+1} \ .
\label{eq:u}
\end{equation}
To get (\ref{eq:u}) we can apply Taylor formula, at $x_{n+1/2}$, to both $u_{n+1}$ and $u_{n}$.
A simpler way to obtain (\ref{eq:u}) is to consider the method of undefined coefficients, for $u_{n+1/2}$ as a linear combination of $u_{n}$ and $u_{n+1}$, and to require that the formula is exact for constant and linear functions.  
In this way we end up with a linear system of two equations in two unknowns where the coefficient matrix is a Vandermonde matrix. 
As far as the first derivative is concerned we can apply the following approximation
\begin{equation}
\frac{du}{dx}(x_{n+1/2}) \approx \frac{u_{n+1}-u_n}{2\left(x_{n+3/4} - x_{n+1/4}\right)} \ .
\label{eq:du}
\end{equation}
These formulae use the value $ u_N = u_\infty $, but not $ x_N = \infty $.
In order to justify non-standard finite difference formula (\ref{eq:du}) we note that, by considering $u=u(\xi(x))$, we can write
\begin{equation}
\left. \frac{du}{dx}\right|_{n+1/2} = \left. \frac{du}{d\xi}\right|_{n+1/2}\left. \frac{d\xi}{dx}\right|_{n+1/2} \approx \frac{u_{n+1}-u_n}{\xi_{n+1} - \xi_{n}} 
\frac{2\left(\xi_{n+3/4}-\xi_{n+1/4}\right)}{2\left(x_{n+3/4} - x_{n+1/4}\right)} \ .
\label{eq:duder}
\end{equation}
The last formula on the right hand side of equation (\ref{eq:duder}) reduces to the right hand side of equation (\ref{eq:du}) because we are using a uniform grid for $\xi$ and therefore $2\left(\xi_{n+3/4}-\xi_{n+1/4}\right)=\xi_{n+1} - \xi_{n}$.

The two finite difference approximations (\ref{eq:u}) and (\ref{eq:du}) have order of accuracy $O(N^{-2})$.
For a system of differential equations, formulae (\ref{eq:u}) and (\ref{eq:du}) can be applied component-wise.

\section{BVPs on infinite intervals}
Let us consider the class of BVPs defined on an infinite interval
\begin{eqnarray}
&& {\ds \frac{d{\bf u}}{dx}} = {\bf f} \left(x, {\bf u}\right)
\ , \quad x \in [0, \infty) \ , \nonumber \\[-1.5ex]
\label{p} \\[-1.5ex]
&& {\bf g} \left( {\bf u}(0), {\bf u} (\infty) \right) = {\bf 0}
\ ,  \nonumber
\end{eqnarray}
where $ {\bf u}(x) $ is a $ d-$dimensional vector with $ u_{\ell}
(x) $ for $ \ell =1, \dots , d $ as components, $ {\bf f}:[0,
\infty) \times \RR^d \rightarrow~\RR^d $, and $ {\bf g}:
\RR^d \times \RR^d \rightarrow \RR^d $.

A non-standard finite difference scheme on a quasi-uniform grid for the class of BVPs
(\ref{p}) can be defined by using the approximations given by (\ref{eq:u}) and (\ref{eq:du}).
 We denote by the $d-$dimensional  vector $ {\bf U}_n $ the numerical approximation to the solution $ {\bf u} (x_n) $ of (\ref{p}) at the points of the mesh, that is for $ n = 0, 1, \dots , N $ .
A second order finite difference scheme for
(\ref{p}) can be written as follows:
\begin{eqnarray}
& {\bf U}_{n+1} - {\bf U}_{n} - a_{n+1/2} {\bf f} \left( x_{n+
1/2}, b_{n+1/2}{\bf U}_{n+1} + c_{n+1/2}{\bf U}_{n} \right) = {\bf 0}
\ , \nonumber\\
& \mbox{for} \quad n=0, 1, \dots , N-1
\label{boxs} \\ 
& {\bf g} \left( {\bf U}_0,{\bf U}_N \right) = {\bf 0} \ ,  \nonumber
\end{eqnarray}
where 
\begin{eqnarray}\label{eq:abc}
a_{n+1/2} &=& 2\left(x_{n+3/4} - x_{n+1/4}\right) \ , \nonumber \\
b_{n+1/2} &=& \frac{x_{n+1/2}-x_{n}}{x_{n+1}-x_n} \ , \\
c_{n+1/2} &=& \frac{x_{n+1}-x_{n+1/2}}{x_{n+1}-x_n} \nonumber \ , 
\end{eqnarray}
for $n=0, 1, \dots , N-1$.

It is evident that (\ref{boxs}) is a nonlinear system of $ d \cdot (N+1)$ equations in the $ d \cdot (N+1)$ unknowns $ {\bf U} = ({\bf U}_0,
{\bf U}_1, \dots , {\bf U}_N)^T $. 
We notice that $b_{n+1/2} \approx c_{n+1/2} \approx 1/2$ for all $n=0, 1, \dots , N-2$, but when $n=N-1$, then $b_{N-1/2} = 0$ and $c_{N-1/2} = 1$. 
This means that $u_{n-1/2} = u_{n-1}$ and the value of $u_N$ plays no role in defining the middle node value.
In order to avoid a sudden jump for the coefficients of (\ref{boxs}), and to make use of $u_N$, we choose to set $b_{N-1/2} = b_{N-3/2}$ and $c_{N-1/2} = c_{N-3/2}$ .
As it will be clear from the results reported in the next section this choice produces a much smaller error in the numerical solution of the system at $x_N$.

For the solution of (\ref{boxs}) we can apply the classical Newton's method along with the simple termination criterion
\begin{eqnarray*}
{\ds \frac{1}{d(N+1)} \sum_{\ell =1}^{d} \sum_{n=0}^{N}
|\Delta U_{n \ell}| \leq {\rm TOL}} \ ,
\end{eqnarray*}
where $ \Delta U_{n \ell} $, for $ n = 0,1, \dots, N $ and $ \ell
= 1, 2, \dots , d $, is the difference between two successive
iterate components and $ {\rm TOL} $ is a fixed tolerance. The
results listed in the next sections were computed by
setting $ {\rm TOL} = 1\E-6 $.

%\section{Test problems}\label{S:test} 
\section{The Falkner-Skan model}
The Falkner-Skan model \cite{Falkner:1931:SAS} of boundary layer theory 
%, see for instance Schlichting 
\cite{Schlichting:2000:BLT}, defined by
\begin{eqnarray}
& {\displaystyle \frac{d^3 u}{d x^3}} + u {\displaystyle \frac{d^{2}u}{dx^2}} + 
P \left[ 1 - \left( {\displaystyle \frac{du}{dx}}\right)^2 \right] = 0 \nonumber \\ [-1ex]
\label{eq:Falkner} \\[-1.5ex]
& u(0) = {\displaystyle \frac{du}{dx}}(0) = 0, \qquad
{\displaystyle \frac{du}{dx}}(\infty) = 1  \ , \nonumber
\end{eqnarray}
is a BVP defined on a semi-infinite interval.

As a first step we rewrite the Falkner-Skan equation in (\ref{eq:Falkner}) as a first order system by setting
\begin{eqnarray*}
u_{i+1}(x) = \ds\frac{d^{i} u}{dx^{i}} (x) \ , \quad
\mbox{for} \ i = 0, 1, 2 \ .
\end{eqnarray*}
In this way the original BVP (\ref{eq:Falkner}) specializes to
\begin{eqnarray}\label{eq:Falkner:system}
&& {\displaystyle \frac{du_1}{dx}} = u_2 \nonumber \\[0.5ex]
&& {\displaystyle \frac{du_2}{dx}} = u_3 \nonumber \\[-1ex]
&& \\[-1ex]
&& {\displaystyle \frac{du_3}{dx}} = -u_1 u_3 - P (1 - {u_2}^2) \nonumber \\[0.5ex]
&&  u_1(0) = u_2(0) = 0 \ ,
\qquad u_2(\infty) = 1 \ , \nonumber
\end{eqnarray}
that is,
\begin{eqnarray*}
& {\bf u} = (u_1,u_2,u_3)^T \\
& {\bf f} (x, {\bf u}) =
\left(u_2,u_3,u_4, -u_1 u_3 - P (1 - {u_2}^2) \right)^T \\
& {\bf g} \left( {\bf u}(0), {\bf u} (\infty) \right) =
\left(u_1(0),u_2(0),u_2(\infty) -1\right)^T
\end{eqnarray*}
in (\ref{p}).

\subsection{Numerical results}
\begin{figure}[!hbt]
\centering
\psfrag{x}[][]{$x$} 
\psfrag{y}[][]{}
%\psfrag{y}[][]{$u(x), {\displaystyle \frac{du}{dx}(x)},{\displaystyle \frac{d^2u}{dx^2}(x)}$} 
\includegraphics[width=.9\textwidth]{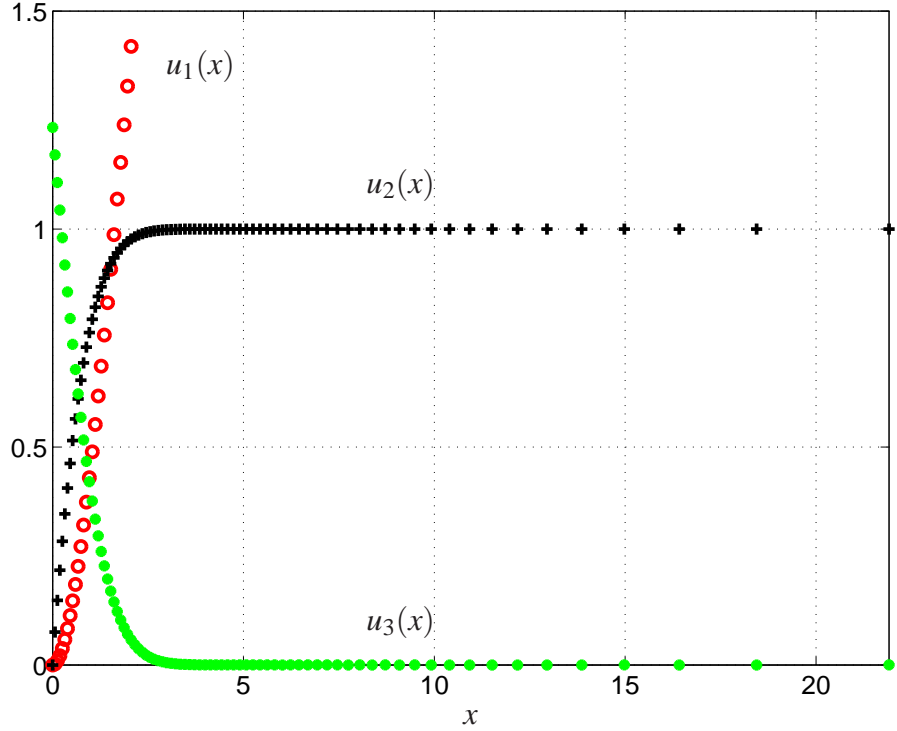}
\put(-275,250){${\displaystyle u_1(x)}$}
\put(-200,205){${\displaystyle u_2(x)}$}
\put(-200,40){${\displaystyle u_3(x)}$}
\caption{\it Numerical solution of Falkner-Skan model with $P=1$ obtained with the map $x=x(\xi)$ defined by (\ref{eq:qu1}) with $c=5$ for $N=80$. }
\label{fig:Falkner}
\end{figure}

%Some of the numerical results obtained for the test problem introduced in this section are reported in this subsection. 
Table \ref{tab:Falkner} lists the numerical approximations of $\displaystyle{\frac{d^2u}{dx^2}}(\infty)$, $\displaystyle{\frac{d^2u}{dx^2}(0)}$ and order of accuracy for increasing values of $N$. 
Here, and in the following, the order of accuracy $p$ is defined by
\begin{equation}\label{eq:p}
p = {\ds \frac{\log(|T_N-T_{1280}|)-\log(|T_{2 N}-T_{1280}|)}{\log(2)}} \ ,
\end{equation}
where $T_N$ and $T_{2 N}$ are numerical approximations of our missing initial condition, namely $\displaystyle{\frac{d^2u}{dx^2}}(0)$.
For all values of $N$ we used the initial iterate
\begin{equation}
u_1(x) = u_2(x) = 1/2 \ , \quad u_3(x) = 10^{-2} \ .
\end{equation}
Let us remark here that the same values of $\displaystyle{\frac{d^2u}{dx^2}}(0)$ were obtained by using $b_{N-1/2} = 0$ and $c_{N-1/2} = 1$ in (\ref{boxs}),
but, on the contrary, larger values of $\displaystyle{\frac{d^2u}{dx^2}}(\infty)$, of the order $O(10^{-3})$, were computed.
In table \ref{tab:Falkner:comp} we compare the obtained numerical results with those available in literature: the agreement is really good.
\begin{table}[!hbt]
\caption{\it Numerical approximation of $\displaystyle{\frac{d^2u}{dx^2}}(\infty)$, $\displaystyle{\frac{d^2u}{dx^2}}(0)$ and order of accuracy.}
\begin{center}
{\begin{tabular}{rcr@{.}lcc}
\hline \\[-1.5ex]
$N$ & iter  & \multicolumn{2}{c}% 
{$\displaystyle{\frac{d^2u}{dx^2}(\infty)}$} 
& $\displaystyle{\frac{d^2u}{dx^2}(0)}$ & {$p$} \\[1.5ex]
\hline
%10    & 5 & $ 0$ & $23 \cdot 10^{-3}$ & $1.256996$   \\  
20    & 6 & $-0$ & $21 \cdot 10^{-7}$ & $1.238724$ &   \\ 
40    & 5 & $ 0$ & $24 \cdot 10^{-7}$ & $1.234124$ & 1.998825 \\
80    & 5 & $-0$ & $33 \cdot 10^{-7}$ & $1.232972$ & 2.002822  \\ 
160   & 5 & $ 0$ & $14 \cdot 10^{-7}$ & $1.232684$ & 2.011345  \\  
320   & 5 & $-0$ & $25 \cdot 10^{-7}$ & $1.232612$ & 2.046294  \\   
640   & 5 & $ 0$ & $39 \cdot 10^{-7}$ & $1.232594$ & 2.201634 \\  
1280  & 5 & $ 0$ & $33 \cdot 10^{-7}$ & $1.232589$ & \\  
\hline
\end {tabular}}
\end{center}
\label{tab:Falkner} 
\end{table}
Figure \ref{fig:Falkner} shows the numerical solution of Falkner-Skan model (\ref{eq:Falkner}) with $P =1$ and $N = 80$.
\begin{table}[!htb]
\caption{\it Comparison of $\displaystyle{\frac{d^2 u}{dx^2}}(0)$
and free or truncated boundary ($ x_\epsilon $ and $x_\infty $ respectively) for the Homann ($P=1/2$) and Hiemenz ($P=1$) flows.}
\begin{center}{ \renewcommand\arraystretch{1.3}
\begin{tabular}{l|r@{.}lccccccc} 
\hline%
& \multicolumn{3}{c}%
{Nasr et al. \cite{Nasr:1990:CSL}}
& \multicolumn{2}{c}%
{Fazio \cite{Fazio:1994:FSEb}} 
& \multicolumn{2}{c}%
{Asaithambi \cite{Asaithambi:1998:FDM}} 
& \multicolumn{2}{c}%
{This paper}\\ 
&  \multicolumn{3}{c}%
{Chebyshev method} & \multicolumn{2}{c}%
{Free BF} & \multicolumn{2}{c}%
{Finite difference} 
& \multicolumn{2}{c}%
{Quasi-uniform} \\ 
\hline\\[-2.5ex]%
{$P$} & \multicolumn{2}{c}%
{$\ \ \ x_\infty $} & 
{$ {\displaystyle \frac{d^2 u}{dx^2}} (0) $} & 
{$ x_\epsilon $} &
{$ {\displaystyle \frac{d^2 u}{dx^2}} (0) $} & 
{$ x_\infty $} &
{$ {\displaystyle \frac{d^2 u}{dx^2}} (0) $} 
& {$x_N$} & {${\displaystyle \frac{d^2 u}{dx^2}} (0) $} \\[1.5ex] \hline
0.5 & \ \ \ 3 & 7 & 0.927805 &        &           &  &  & & \\
0.5 & \ \ \ 7 & 4 & 0.927680 & 5.09 & 0.927680  & 5.67 & 0.927682 & $\infty$ & 0.927681 \\
1   &  \ \ \ 3 & 5 & 1.232617 &     &           &  &  & & \\
1  &  \ \ \ 7 &   & 1.232588 &  5.19 & 1.232588  & 5.14 & 1.232589  & $\infty$ & 1.232589  \\
\hline
\end{tabular}}
\label{tab:Falkner:comp}
\end{center}
\end{table}

The Falkner-Skan model (\ref{eq:Falkner})  was solved by Grosch and Orszag \cite{Grosch:NSP:1977}, but only for small values of $P$ (namely $-0.1$, $0$, and $0.1$).
Moreover, these authors used the truncated boundary and the two maps (\ref{eq:qu1}) and (\ref{eq:qu2}) with a shooting method using only $11$ mesh point but, unfortunately, reported only the values of the first derivative at $x=1$ instead of the missing initial condition.
Moreover, for $P=-0.1$ they missed the inverse flow solution of Stewartson \cite{Stewartson:1954:FSF,Stewartson:1964:TLB}, see also Asaithambi \cite{Asaithambi:1997:NMS}, Auteri et al. \cite{Auteri:2012:GLS} and Fazio \cite{Fazio:2013:BPF}. 

\section{A pile in soil}
Here we consider a problem that was already used by
Lentini and Keller \cite{Lentini:BVP:1980} 
to test the asymptotic boundary
conditions approach.
This problem is of special interest here because none of
the solution components is a  monotone function,
on $ [0, \infty) $, see \cite{Fazio:2002:SFB,Fazio:2003:FBA}. 
Let $ u(x) $ be the deflection of a semi-infinite pile
embedded in soft soil at a distance $ x $ below the surface of the
soil. The governing differential equation for the movement of the
pile, in dimensionless form, is given by:
\begin{eqnarray}
{\ds \frac{d^4 u}{dx^4}} = - P_1 \left(1 - e^{-P_2 u} \right) \ ,
\qquad x \in \left[0,\infty\right) \ , \nonumber \\[-1.5ex]
\label{eq:pile} \\[-1.5ex]
{\ds \frac{d^2 u}{dx^2}} (0) = 0 \ , \qquad {\ds
\frac{d^3 u}{dx^3}} (0) = P_3 \ , \qquad u(\infty) = {\ds\frac{du}{dx} (\infty)} = 0 \ , \nonumber
\end{eqnarray}
where $ P_1 $ and $ P_2 $ are positive material constants. 
As far as the boundary conditions are concerned, at the
origin a zero moment and a positive shear $P_3$ are assumed and 
from physical considerations it follows that $ u(x) $
and all its derivatives go to zero at infinity,
so that, the zero asymptotic boundary conditions can be imposed. 
This problem is of interest in foundation
engineering: for instance, in the design of drilling rigs above
the ocean floor.

The governing differential equation in (\ref{eq:pile}) can be rewritten as a first order system by setting:
\begin{eqnarray*}
u_{i+1}(x) = \ds\frac{d^{i} u}{dx^{i}} (x) \ , \quad
\mbox{for} \ i = 0, 1, 2, 3 \ .
\end{eqnarray*}
In this way the original BVP (\ref{eq:pile}) specializes to
\begin{eqnarray}
&& {\ds \frac{du_1}{dx}} = u_2 \  , \nonumber \\
&& {\ds \frac{du_2}{dx}} = u_3 \  , \nonumber \\ && {\ds
\frac{du_3}{dx}} = u_4 \ , \label{bvp2}  \\ && {\ds
\frac{du_4}{dx}} = -P_1 \left(1 - {\ds e^{-P_2 u_1}} \right) \  ,
\nonumber \\[0.5ex]
&&  u_3(0) = 0 \ , \qquad u_4(0) = P_3 \ , \qquad u_1(\infty) = 0 \ ,
\qquad u_2(\infty) = 0 \ , \nonumber
\end{eqnarray}
that is,
\begin{eqnarray*}
& {\bf u} = (u_1,u_2,u_3,u_4)^T \\
& {\bf f} (x, {\bf u}) =
\left(u_2,u_3,u_4, -P_1
\left(1 - e^{-P_2 u_1} \right)\right)^T \\
& {\bf g} \left( {\bf u}(0), {\bf u} (\infty) \right) =
\left(u_3(0),u_4(0) - P_3,u_1(\infty),u_2(\infty) \right)^T
\end{eqnarray*}
in (\ref{p}).

\subsection{Numerical results}
%In this subsection we report some of the numerical results obtained for the test problem introduced in this section. 
In order to be able to compare our numerical results we used the same parameter values employed by Lentini and Keller \cite{Lentini:BVP:1980}
\begin{eqnarray*}
P_1 = 1, \quad P_2 = \frac{1}{2} \quad \mbox{and}
\quad P_3 = \frac{1}{2} \ .
\end{eqnarray*}
Moreover, we choose to consider the values of the missing initial
conditions $ u_1 (0) $ and $ u_2 (0) $ as representative results.

\begin{figure}[!hbt]
\centering
\psfrag{x}[][]{$x$} 
\psfrag{y}[][]{}
%\psfrag{y}[][]{$u(x), {\displaystyle \frac{du}{dx}(x)},{\displaystyle \frac{d^2u}{dx^2}(x)},{\displaystyle \frac{d^3u}{dx^3}(x)}$} 
\includegraphics[width=.9\textwidth]{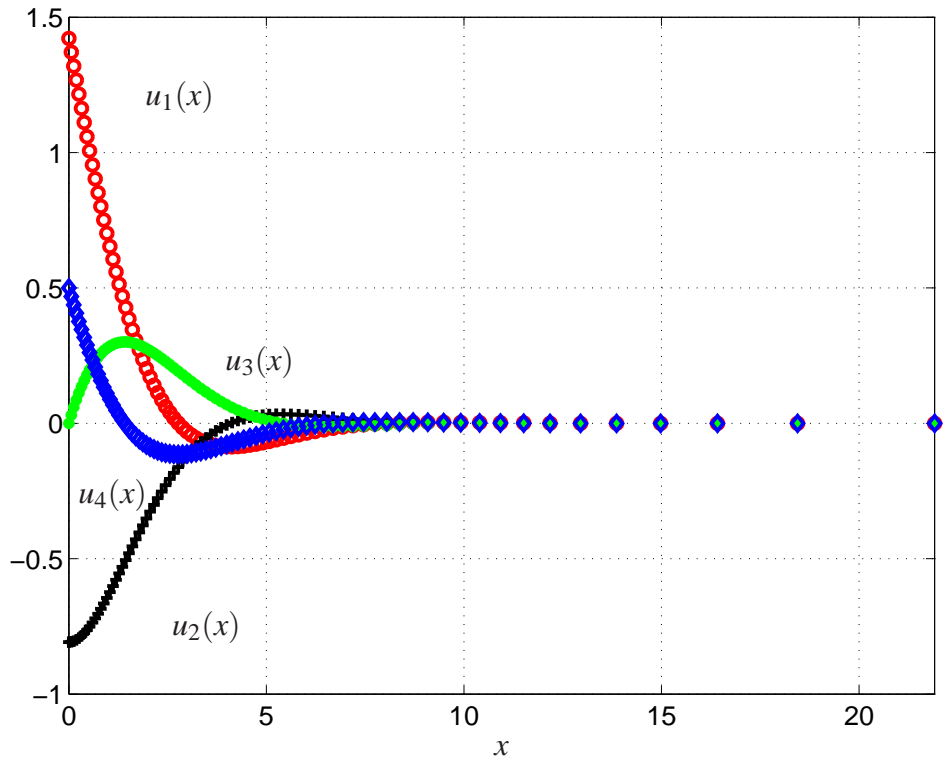}
\put(-300,250){${\displaystyle u_1(x)}$}
\put(-290,50){${\displaystyle u_2(x)}$}
\put(-270,150){${\displaystyle u_3(x)}$}
\put(-325,100){${\displaystyle u_4(x)}$}
\caption{\it Numerical solution of pile model (\ref{eq:pile}) obtained with the map $x=x(\xi)$ defined by (\ref{eq:qu1}) with $c=5$ for $N=80$. }
\label{fig:pile}
\end{figure}

Figure \ref{fig:pile} shows the numerical solution of the BVP (\ref{eq:pile}) using the map (\ref{eq:qu1}) with $c=5$ for $N=80$. %obtained with $P_1=1$, $P_2 = 0.5$ and $P_3 = 0.5$.

Table \ref{tab:pile} lists the numerical approximations of $u(0)$, $\displaystyle{\frac{du}{dx}(0)}$ and the corresponding order of accuracy for increasing values of $N$. 
The numerical accuracy values $p$, once again, are computed by formula (\ref{eq:p}),
where $T_N$ and $T_{2 N}$ are numerical approximations of our missing initial conditions, namely $u(0)$ and $\displaystyle{\frac{du}{dx}}(0)$, respectively.
For all values of $N$ we used the initial iterate
\[
u_1(x) = u_2(x) = u_3(x) = u_4(x) = 1 \ .
\]

\begin{table}[!hbt]
\caption{\it Numerical approximation of $u(0)$, $\displaystyle{\frac{du}{dx}}(0)$ and the corresponding order of accuracy. {\rm NaN} means not a number.}
\begin{center}
{\begin{tabular}{rccccc}
\hline \\[-1.5ex]
$N$ & iter & $u(0)$ & {$p$} & $-\displaystyle{\frac{du}{dx}(0)}$ & {$p$} \\[1.5ex]
\hline
%10    & 5 & $1.416457$ & & $0.804538$ & \\  
20    & 5 & $1.420337$ &            & $0.807289$ & \\ 
40    & 5 & $1.421243$ & $2.003590$ &$0.807934$ & $2.020368$ \\
80    & 5 & $1.421469$ & $2.004801$ & $0.808094$ & $2.048674$ \\  
160   & 5 & $1.421526$ & $2.058894$ & $0.808135$ & $2.350497$ \\ 
320   & 5 & $1.421540$ & $2.169925$ & $0.808145$ & {$\infty$} \\ 
640   & 5 & $1.421544$ & {$\infty$} & $0.808145$ & {\rm NaN} \\ 
1280  & 5 & $1.421544$ &            & $0.808145$ & \\ 
\hline
\end {tabular}}
\end{center}
\label{tab:pile} 
\end{table}

\begin{table}[!htb]
\caption{\it Comparison of $u(0)$, $\displaystyle{\frac{du}{dx}}(0)$
and free or truncated boundary ($ x_\epsilon $ and $x_\infty $ respectively) for the pile problem.}
\begin{center}{ \renewcommand\arraystretch{1.3}
\begin{tabular}{cccccc} 
\hline%
\multicolumn{2}{c}%
{Lentini and Keller \cite{Lentini:BVP:1980}} 
&\multicolumn{2}{c}%
{Fazio \cite{Fazio:1994:FSEb}} 
& \multicolumn{2}{c}%
{This paper} \\ 
 \multicolumn{2}{c}%
{Asymptotic BCs {$ x_\infty = 10$}} 
& \multicolumn{2}{c}%
{Free BF {$ x_\epsilon = 17.75$}} 
&\multicolumn{2}{c}%
{Quasi-uniform {$x_N=\infty$}}
 \\ 
\hline\\[-2.5ex]%
{$u(0) $} & {$ {\displaystyle \frac{du}{dx}}(0) $} 
& {$u(0) $}
& {$ {\displaystyle \frac{du}{dx}}(0) $}  
& {$u(0) $}
& {$ {\displaystyle \frac{du}{dx}}(0) $}
 \\[1.5ex] \hline
$1.4215$ & $-0.80814$ & 1.42154  & $ -0.808144 $ & $1.421544$ & $-0.808145$  \\
\hline
\end{tabular}}
\label{tab:pile:comp}
\end{center}
\end{table}
As far as the missing initial conditions for the pile problem are
concerned, a comparison of the obtained values has been considered in table \ref{tab:pile:comp}. 
It is easily seen that our results are in good agreement with those available in literature.
%of $ u_1 (0) $ and $ u_2 (0) $ can be made, respectively, with the values $ 1.4215 $ and $ -.80814 $ reported by 
Lentini and Keller \cite{Lentini:BVP:1980} used the mentioned asymptotic boundary conditions and employed PASVAR, a routine based upon the
trapezoidal difference scheme with automatic mesh refinement and
deferred corrections as described by Lentini and Pereyra \cite{Lentini:AFD:1977}. 
That software is rather sophisticated because it
adjusts automatically the mesh and the order of accuracy of the method employed.
Fazio \cite{Fazio:1994:FSEb} used a free boundary formulation of the pile problem as mentioned in the introduction and the Keller's box finite difference method.

\subsection{Improving the accuracy via Richardson's extrapolation}
There are two possible strategies for improving the accuracy of the obtained numerical results.
The first one is to study higher order finite difference schemes on quasi-uniform grids, whereas the second one is to apply Richardson's extrapolation.
As we shall see shortly the simplest of these two strategies is to apply Richardson's extrapolation, and we will explain in full details its application to the numerical results obtained by our scheme on a quasi-uniform grid. 

However, let us indicate first the steps necessary to develop higher order schemes.
First of all, we have to define the values of $u(x)$ on the mid-points of the grid using a wider stencil
\begin{equation}
u_{n+1/2} \approx \alpha_{-1} u_{n-1} + \alpha_0 u_n + \alpha_1 u_{n+1} +\alpha_2 u_{n+2} \ ,
\label{eq:u4ord}
\end{equation}
where $\alpha_j$ for $j = -1, 0, 1, 2$ are constants to be determined.
To get the coefficients in (\ref{eq:u4ord}) we require that the formula is exact for constant, linear, quadratic and cubic functions.  
In this way we end up with the linear system of four equations in four unknowns $V\alpha = d$, where
\begin{equation}\label{eq:matrix}
V =
\begin{bmatrix}
1 & 1 & 1 & 1 \\
x_{n-1} & x_{n} & x_{n+1} & x_{n+2} \\
x^2_{n-1} & x^2_{n} & x^2_{n+1} & x^2_{n+2} \\
x^3_{n-1} & x^3_{n} & x^3_{n+1} & x^3_{n+2} 
\end{bmatrix}
\end{equation}
$\alpha = [\alpha_{-1},\alpha_0, \alpha_{1}, \alpha_2]^T$ and $d = [1,x_{n+1/2}, x^2_{n+1/2}, x^3_{n+1/2}]^T$.
We notice that the grid-points are all distinct, hence the coefficient matrix is a Vandermonde matrix and the solution of this system exists and is unique. 
However, Vandermonde systems might be ill-conditioned, see Gautschi \cite{Gautschi:IVC:1962,Gautschi:NEI:1975} and his review paper \cite{Gautschi:HUA:1990}.
In fact, the condition number of a Vandermonde matrix may be large, causing large errors when computing the solution of the system numerically.

As far as the higher order finite difference formula for the first derivative is concerned, we have to set, again, a wider stencil
\begin{equation}
\frac{du}{dx}(x_{n+1/2}) \approx \beta_{-1} u_{n-1} + \beta_0 u_n + \beta_1 u_{n+1} +\beta_2 u_{n+2} \ ,
\label{eq:du4ord}
\end{equation}
where $\beta_j$ for $j = -1, 0, 1, 2$ are constants to be determined.
Once again, in order to find the coefficients in (\ref{eq:du4ord}) we require that the formula is exact for constant, linear, quadratic and cubic functions.  
In this way we end up with the linear system of four equations in four unknowns $V\beta = r$, where $V$ is the same matrix defined above in (\ref{eq:matrix}),
$\beta = [\beta_{-1},\beta_0, \beta_{1}, \beta_2]^T$ and $r = [0,1,2 x_{n+1/2}, 3 x^2_{n+1/2}]^T$.
Of course, both systems can be solved exactly, for instance by using a general purpose Computer Algebra system, like the free software AXIOM \cite{AXIOM} or Derive \cite{Derive}. 
However, the application of higher order finite difference schemes is not straightforward because we have to get and implement also boundary finite difference formulas of the same order.
Since the computational stencil is wider than in the second order case we need to introduce ghost cells at the boundary and this means that we need a ghost cell greater than infinity.
If we apply lower order boundary conditions, then a reduction of the overall accuracy results, see Fazio and Russo \cite{Fazio:2010:LCS} for numerical results related to the numerical implementation of higher order boundary conditions in the study the dynamics of two gases in a piston problem. 

At this stage we prefer to indicate a different strategy to get higher accuracy.
In fact, this is a further advantage in using a family of quasi-uniform grids in calculations. 
The algorithm is based on Richardson's extrapolation, introduced by Richardson in \cite{Richardson:1927:DAL}, and it is the same for many
finite difference methods: for numerical differentiation or integration, solving systems of ordinary or partial differential equations.
To apply Richardson's extrapolation, we carry on two calculations on embedded uniform or quasi-uniform grids with total number of nodes $N$ and $2N$. 
All nodes of largest steps are identical to even nodes of denser grid due to uniformity. 
Let us suppose that we use a numerical method with order of accuracy $p$ to find an approximation of a scalar value $T$. 
For smooth enough solutions the error can be decomposed into a sum of inverse powers of $N$. 
The Richardson's formula
\begin{equation}
R_{2N} = \frac{T_{2N}-T_N}{2^p-1} 
\label{eq:Retra}
\end{equation}
defines the main term of such sum. 
This formula is asymptotically exact in the limit as $N$ goes to infinity if we use uniform or quasi-uniform grids. 
Hence, it gives the real value of numerical solution error without knowledge of exact solution.
We can apply (\ref{eq:Retra}) as a single-step correction formula
\begin{equation}
T = T_{2N} + R_{2N} + O(N^{-p-s})
\label{eq:conextra}
\end{equation}
and increase the order of accuracy of our approximation. 
In general, we have $s = 2$ in (\ref{eq:conextra}) for a symmetrical finite difference scheme, but only $s = 1$ for non-symmetrical ones. 
Such enlargement of accuracy requires to solve the same problem twice and only few further arithmetical operations and so it is very cheap.

It could be possible to apply Richardson's extrapolation to the results reported in table \ref{tab:Falkner} and in table \ref{tab:pile}.
This has been done in table \ref{tab:Falkner:extra} and in tables \ref{tab:pile:extra1} and \ref{tab:pile:extra2}, respectively.
\begin{table}[!hbt]
\caption{\it Richardson's extrapolation for the first line of table \ref{tab:Falkner}.}
\begin{center}
{\begin{tabular}{rcc}
\hline \\[-1.5ex]
$N$ & $T^{(0)}$ & $T^{(1)}$  \\[1.5ex]
\hline
%20    & $1.238724$ & $ $ \\
40    & $1.234124$ &  \\
80    & $1.232972$ & $ 1.232588 $ \\ 
160   & $1.232684$ & $ 1.232588 $ \\  
\hline
\end {tabular}}
\end{center}
\label{tab:Falkner:extra} 
\end{table}

\begin{table}[!hbt]
\caption{\it Richardson's extrapolation of $u(0)$ for the first line of table \ref{tab:pile}.}
\begin{center}
{\begin{tabular}{rccc}
\hline \\[-1.5ex]
$N$ & $T^{(0)}$ & $T^{(1)}$ & $T^{(2)}$ \\[1.5ex]
\hline
%20    & $1.420337$ & $ $ & \\ 
40    & $1.421243$ &  & \\
80    & $1.421469$ & $1.421544$ & \\ 
160   & $1.421526$ & $1.421545$ & $1.421545$ \\  
\hline
\end {tabular}}
\end{center}
\label{tab:pile:extra1} 
\end{table}

\begin{table}[!hbt]
\caption{\it Richardson's extrapolation of $\displaystyle{\frac{du}{dx}}(0)$ for the first line of table \ref{tab:pile}.}
\begin{center}
{\begin{tabular}{rccc}
\hline \\[-1.5ex]
$N$ & $T^{(0)}$ & $T^{(1)}$ & $T^{(2)}$ \\[1.5ex]
\hline
40    & $-0.807934$ & $ $ & \\
80    & $-0.808094$ & $-0.808147$ & \\ 
160   & $-0.808135$ & $-0.808149$ & $-0.808149$\\  
\hline
\end {tabular}}
\end{center}
\label{tab:pile:extra2} 
\end{table}

In this way we can avoid to solve the given problem with a large number of grid-points.

A simple extrapolation improves the numerical accuracy obtained for the results to a given problem.
Of course, it is also possible to iterate the extrapolation until $T^{(k)}_{2N} = T^{(k)}_{N}$, for $k = 1,2,\dots$, as in table \ref{tab:Falkner:extra}.
We also stop to apply a nested extrapolation as soon as $T^{(k)}_{N} = T^{(k-1)}_N$ as in tables \ref{tab:pile:extra1} and \ref{tab:pile:extra2}.

Starting with the computed data $T^{(0)}_{}$, we proceed with the extrapolated values
\begin{equation}
T^{(k)} = \frac{2^{k+1} T^{(k-1)}_{2N}-T^{(k-1)}_N}{2^{k+1}-1}  \ ,
\label{eq:stpRetra1}
\end{equation}
for $k = 1, 2, \dots $.
The extrapolated values reported in tables \ref{tab:Falkner:extra}, \ref{tab:pile:extra1} and \ref{tab:pile:extra2} can be compared with the corresponding values
listed in tables \ref{tab:Falkner} and \ref{tab:pile}, respectively.
It is clear that to get accurate numerical results, for both the considered problems, we can apply few extrapolations with nested quasi-uniform grids involving a small number of grid points.
%and the second by
%\begin{equation}
%T^{(2)} = \frac{8 T^{(1)}_{2N}-T^{(1)}_N}{7} \ .
%\label{eq:stpRetra2}
%\end{equation}
%Of course, it is an easy matter to define higher order extrapolation formulae, by
%\begin{equation}
%T^{(3)} = \frac{16 T^{(2)}_{2N}-T^{(2)}_N}{15} \ .
%\label{eq:stpRetra3}
%\end{equation}

\section{Concluding remarks}
The numerical results for the test problems reported in the previous sections show that non-standard finite difference schemes on quasi-uniform grids are an effective way to solve BVPs defined on infinite intervals.
The application of non-standard finite difference schemes on quasi-uniform grids overcomes the need for a priori definition of the truncated boundary. 
For both examples we used the logarithmic map (\ref{eq:qu1}) because in both cases the largest variation of the solution components occur near the origin.

Let us discuss, at the end of this work, a possible way to extend
the non-standard finite difference schemes on quasi-uniform grids to the numerical solutions of problems
defined on the whole real line, for instance, the connecting orbits
problems mentioned in the introduction.
For these problems, all boundary conditions are imposed at plus or minus infinity.
In such a case it is possible to use the tangential quasi-uniform grid 
\[
x_n = c \cdot \tan\left(\frac{n \pi}{2N}\right) \ ,
\] 
where $c > 0$ is a control parameter.
In fact, if $ n = -N,-N+1,\dots -1, 0, 1, \dots , N-1, N$, then this tangential grid cover the whole infinite line, and in particular we have that $x_{-N} = -\infty$.

\vspace{1.5cm}

\noindent {\bf Acknowledgement.} {The research of this work was 
supported, in part, by the University of Messina and by the GNCS of INDAM.}


\begin{thebibliography}{10}

\bibitem{Asaithambi:1997:NMS}
A.~Asaithambi.
\newblock A numerical method for the solution of the {F}alkner{-}{S}kan
  equation.
\newblock {\em Appl. Math. Comput.}, 81:259--264, 1997.

\bibitem{Asaithambi:1998:FDM}
A.~Asaithambi.
\newblock A finite-difference method for the solution of the {F}alkner{-}{S}kan
  equation.
\newblock {\em Appl. Math. Comput.}, 92:135--141, 1998.

\bibitem{Auteri:2012:GLS}
F.~Auteri and L.~Quartapelle.
\newblock Galerkin-laguerre spectral solution of self-similar boundary layer
  problems.
\newblock {\em Commun. Comput. Phys.}, 12:1329--1358, 2012.

\bibitem{AXIOM}
AXIOM home-page: \url{http://www.axiom-developer.org/}

\bibitem{Bai:1993:NCH}
F.~Bai, A.~Spence, and A.~M. Stuart.
\newblock The numerical computation of heteroclinic connections in systems of
  gradient partial differential equations.
\newblock {\em SIAM J. Appl. Math.}, 53:743--769, 1993.

\bibitem{Beyn:1990:GBN}
W.~J. Beyn.
\newblock Global bifurcation and their numerical computation.
\newblock In D.~Rossed, B.~D. Dier, and A.~Spence, editors, {\em Bifurcation:
  Numerical Techniques and Applications}, pages 169--181. Kluwer, Dordrecht,
  1990.

\bibitem{Beyn:1990:NCC}
W.~J. Beyn.
\newblock The numerical computation of connecting orbits in dynamical systems.
\newblock {\em IMA J. Numer. Anal.}, 9:379--405, 1990.

\bibitem{Beyn:1992:NMD}
W.~J. Beyn.
\newblock Numerical methods for dynamical systems.
\newblock In W.~Light, editor, {\em Advances in Numerical Analysis}, pages
  175--236. Clarendon Press, Oxford, 1992.

\bibitem{Botta:1972:NSN}
E.~F.~F. Botta, D.~Dijkstra, and A.~E.~P. Veldman.
\newblock The numerical solution of the {N}avier-{S}tokes for laminar,
  incompressible flow past a parabolic cylinder.
\newblock {\em J. Eng. Math.}, 6:63--81, 1972.

\bibitem{Boyd:2001:CFS}
J.~P. Boyd.
\newblock {\em Chebyshev and Fourier Spectral Methods}.
\newblock Dover, New York, 2001.

\bibitem{Collatz}
L.~Collatz.
\newblock {\em The Numerical Treatment of Differential Equations}.
\newblock Springer, Berlin, 3rd edition, 1960.

\bibitem{Davis:1972:NSN}
R.~T. Davis.
\newblock Numerical solution of the {N}avier-{S}tokes equations for symmetric
  laminar incompressible flow past a parabola.
\newblock {\em J. Fluid Mech.}, 51:417--433, 1972.

\bibitem{deHoog:1980:ATB}
F.~R. de~Hoog and R.~Weiss.
\newblock An approximation theory for boundary value problems on infinite
  intervals.
\newblock {\em Computing}, 24:227--239, 1980.

\bibitem{Derive}
Derive home-page: \url{http://derive.en.softonic.com/}

\bibitem{Falkner:1931:SAS}
V.~M. Falkner and S.~W. Skan.
\newblock Some approximate solutions of the boundary layer equations.
\newblock {\em Philos. Mag.}, 12:865--896, 1931.

\bibitem{Fazio:1992:BPF}
R.~Fazio.
\newblock The {Blasius} problem formulated as a free boundary value problem.
\newblock {\em Acta Mech.}, 95:1--7, 1992.

\bibitem{Fazio:1994:FSEb}
R.~Fazio.
\newblock The {Falkner}-{Skan} equation: numerical solutions within group
  invariance theory.
\newblock {\em Calcolo}, 31:115--124, 1994.

\bibitem{Fazio:1996:NAN}
R.~Fazio.
\newblock A novel approach to the numerical solution of boundary value problems
  on infinite intervals.
\newblock {\em SIAM J. Numer. Anal.}, 33:1473--1483, 1996.

\bibitem{Fazio:2002:SFB}
R.~Fazio.
\newblock A survey on free boundary identification of the truncated boundary in
  numerical {BVP}s on infinite intervals.
\newblock {\em J. Comput. Appl. Math.}, 140:331--344, 2002.

\bibitem{Fazio:2003:FBA}
R.~Fazio.
\newblock A free boundary approach and {K}eller's box scheme for {BVP}s on
  infinite intervals.
\newblock {\em Int. J. Computer Math.}, 80:1549--1560, 2003.

\bibitem{Fazio:2013:BPF}
R.~Fazio.
\newblock {B}lasius problem and {F}alkner-{S}kan model: {T}{\"o}pfer's
  algorithm and its extension.
\newblock {\em Comput. \& Fluids}, 73:202--209, 2013.
\newblock Preprint available at
  the URL: \url{http://arxiv.org/pdf/1212.5057v1.pdf}.

\bibitem{Fazio:2014:QUG}
R.~Fazio and A.~Jannelli.
\newblock Quasi-uniform grids and ad hoc finite difference schemes for bvps on
  infinite intervals.
\newblock ENUMATH Conference 2013, Lausanne, August 26-30, 2013.
\newblock Presentation available at the URL: \url{http://mat521.unime.it/~fazio/preprints/ENUMATH2013_pdfscreen.pdf}

\bibitem{Fazio:2010:LCS}
R.~Fazio and G.~Russo.
\newblock Central schemes and second order boundary conditions for 1{D}
  interface and piston problems in lagrangian coordinates.
\newblock {\em Commun. Comput. Phys.}, 8:797--822, 2010.

\bibitem{Fox}
L.~Fox.
\newblock {\em Numerical Solution of Two-point Boundary Value Problems in
  Ordinary Differential Equations}.
\newblock Clarendon Press, Oxford, 1957.

\bibitem{Friedman:1991:NCC}
M.~J. Friedman and E.~J. Doedel.
\newblock Numerical computation and continuation of invariant manifolds
  connecting fixed points.
\newblock {\em SIAM J. Numer. Anal.}, 28:789--808, 1991.

\bibitem{Gautschi:IVC:1962}
W. Gautschi.
\newblock On inverses of Vandermonde and confluent Vandermonde matrices.
\newblock {\em Numer. Math.}, 4:117--123, 1962.

\bibitem{Gautschi:NEI:1975}
W. Gautschi.
\newblock Norm estimates for inverses of Vandermonde matrices.
\newblock {\em Numer. Math.}, 23:337--347, 1975.

\bibitem{Gautschi:HUA:1990}
W. Gautschi.
\newblock How (un)stable are Vandermonde systems?
\newblock Asymptotic And Computational Analysis
(R. Wong, ed.), Lecture Notes in Pure and Applied Mathematics, vol. 124, Marcel Dekker,
Inc., New York and Basel, pp. 193--210. 1990.

\bibitem{Grosch:NSP:1977}
C.~E. Grosch and S.~A. Orszag.
\newblock Numerical solution of problems in unbounded regions: Coordinate
  transforms.
\newblock {\em J. Comput. Phys.}, 25:273--296, 1977.

\bibitem{Koleva:NSH:2006}
M.~N. Koleva.
\newblock Numerical solution of the heat equation in unbounded domains using
  quasi-uniform grids.
\newblock In I.~Lirkov, S.~Margenov, and J.~Wasniewski, editors, {\em
  Large-scale Scientific Computing}, volume 3743 of {\em Lecture Notes in
  Comput. Sci.}, pages 509--517, 2006.

\bibitem{Lentini:BVP:1980}
M.~Lentini and H.~B. Keller.
\newblock Boundary value problems on semi-infinite intervals and their
  numerical solutions.
\newblock {\em SIAM J. Numer. Anal.}, 17:577--604, 1980.

\bibitem{Lentini:AFD:1977}
M.~Lentini and V.~Pereyra.
\newblock An adaptive finite difference solver for two-point boundary value
  problems with mild boundary layers.
\newblock {\em SIAM J. Numer. Anal.}, 14:91--111, 1977.

\bibitem{Liu:1997:CCH}
L.~Liu, G.~Moore, and R.~D. Russell.
\newblock Computation and continuation of homoclinic and heteroclinic orbits
  with arclength parameterization.
\newblock {\em SIAM J. Sci. Comput.}, 18:69--93, 1997.

\bibitem{Markowich:TAS:1982}
P.~A. Markowich.
\newblock A theory for the approximation of solution of boundary value problems
  on infinite intervals.
\newblock {\em SIAM J. Math. Anal.}, 13:484--513, 1982.

\bibitem{Markowich:ABV:1983}
P.~A. Markowich.
\newblock Analysis of boundary value problems on infinite intervals.
\newblock {\em SIAM J. Math. Anal.}, 14:11--37, 1983.

\bibitem{Nasr:1990:CSL}
H.~Nasr, I.~A. Hassanien, and H.~M. El-{H}awary.
\newblock Chebyshev solution of laminar boundary layer flow.
\newblock {\em Int. J. Computer Math.}, 33:127--132, 1990.

\bibitem{Ockendon}
J.~R. Ockendon.
\newblock Math. Rev., 84c:34021.

\bibitem{Richardson:1927:DAL}
L.~F. Richardson.
\newblock The deferred approach to the limit.
\newblock {\em Phil.Trans., A}, 226:299--349, 1927.

\bibitem{Schlichting:2000:BLT}
H.~Schlichting and K.~Gersten.
\newblock {\em Boundary Layer Theory}.
\newblock Springer, Berlin, 8th edition, 2000.

\bibitem{Stewartson:1954:FSF}
K.~Stewartson.
\newblock Further solutions of the {F}alkner-{S}kan equation.
\newblock {\em Proc. Camb. Philos. Soc.}, 50:454--465, 1954.

\bibitem{Stewartson:1964:TLB}
K.~Stewartson.
\newblock {\em The Theory of Laminar Boundary Layers in Compressible Fluids}.
\newblock Oxford University Press, Oxford, 1964.

\bibitem{vandeVooren:1970:NSS}
A.I. van~de Vooren and D.~Dijkstra.
\newblock The {N}avier-{S}tokes solution for laminar flow past a semi-infinite
  flat plate.
\newblock {\em J. Eng. Math.}, 4:9--27, 1970.

\end{thebibliography}
\end{document}